\documentclass[12pt]{article}
\usepackage[dvips]{graphics}
\newcommand{\nc}{\newcommand}
\nc{\dl}{\delta} \nc{\ep}{\varepsilon} \nc{\ph}{\Phi_n}
\renewcommand{\S}{${\cal{S}}_N\>$}
\nc{\bq}{\begin{equation}} \nc{\eq}{\end{equation}} \nc{\vp}{\varphi}
\nc{\ov}{\over} \nc{\hf}{{1\ov2}} \renewcommand{\th}{\theta}
\nc{\inv}{^{-1}}
\renewcommand{\sp}{\vspace{1ex}} \nc{\iy}{\infty} \nc{\cd}{\cdots}
\nc{\La}{\Lambda} \nc{\tn}{\otimes} \nc{\noi}{\noindent}
\nc{\cO}{{\cal{O}}}
\nc{\ba}{\begin{array}} \nc{\ea}{\end{array}} \nc{\pl}{\partial}
\nc{\Tn}{T_n(f)} \nc{\Tni}{T_n(f)\inv} \nc{\dt}{{d\ov dt}}
\nc{\twotwo}[4]{\left(\begin{array}{cc}#1&#2\\&\\#3&#4\end{array}\right)}
\nc{\twoone}[2]{\left(\begin{array}{c}#1\\\\#2\end{array}\right)}
\nc{\ra}{\rightarrow} \nc{\al}{\alpha} \newcommand{\tr}{\mbox{\rm tr}}
\nc{\cU}{{\cal{U}}} \nc{\si}{\sigma} \nc{\E}{E_n} \nc{\la}{\lambda}
\begin{document}

\begin{center} {\large \bf Random Unitary Matrices, Permutations and
Painlev\'e }\end{center}\sp

\sp\begin{center}{{\bf Craig A.~Tracy}\\
{\it Department of Mathematics and Institute of Theoretical Dynamics\\
University of California, Davis, CA 95616, USA\\
e-mail address: tracy@itd.ucdavis.edu}}\end{center}
\begin{center}{{\bf Harold Widom}\\
{\it Department of Mathematics\\
University of California, Santa Cruz, CA 95064, USA\\
e-mail address: widom@math.ucsc.edu}}\end{center}\sp

\sp\begin{center}{\bf Abstract}\end{center}\sp
This paper is concerned with certain connections between the ensemble of
$n\times n$  unitary matrices---specifically the characteristic function
of the
random variable $\tr(U)$---and combinatorics---specifically Ulam's
problem
concerning the  distribution of the length of the longest increasing
subsequence in permutation groups---and the appearance of
Painlev\'e functions in the answers
to apparently  unrelated questions. Among the results is a
representation in
terms of a Painlev\'e V function for the characteristic function  of
$\tr(U)$  and
(using recent results of Baik, Deift and Johansson) an expression in
terms of a Painlev\'e II function for the limiting distributiuon
of the length of the
longest increasing subsequence in the hyperoctahedral groups.
\sp

\renewcommand{\theequation}{1.\arabic{equation}}
\begin{center} {\bf I. Introduction}\end{center}\sp

The characteristic function of the random variable $\tr \,U$, where $U$
belongs to
the ensemble $\cU(n)$ of $n\times n$ unitary matrices with Haar measure,
is the expected value
\bq\E(e^{r\,\tr\,U+s\,\overline{\tr\,U}}).\label{dis}\eq
In $\cU(n)$ we have for any function $g$ with Fourier coefficients $g_k$
\bq\E\left(\prod_{j=1}^n g(e^{i\th_j})\right)=\det\,T_n(g),\label{ev}\eq
where $T_n(g)$ is the associated $n\times n$ Toeplitz matrix defined by
\[T_n(g)=(g_{j-k}),\ \ \ \ \ (j,\,k=0,\cd,n-1).\]
It follows that the distribution function (\ref{dis}) equals the
determinant
of the Toeplitz matrix associated with the function $e^{rz+sz\inv}$.
The determinant, which we denote by $D_n$, is a function of the product
$rs$
(see Sec. II) and so it is completely determined by its values when
$r=s=t$.
This function $D_n(t)$ has connections with both integrable
systems and combinatorial theory.
To state our results, and these connections, we introduce some notation.

We set
\[f(z)=e^{t\,(z+z\inv)},\]
so that $D_n(t)=\det\,\Tn$. Notice that $\Tn$ is symmetric since
$f(z\inv)=f(z)$.
We introduce the $n$-vectors
\[\dl^+=\left(\ba{c}1\\0\\\vdots\\0\\0\ea\right),\ \ \
\dl^-=\left(\ba{c}0\\0\\\vdots\\0\\1\ea\right),
\ \ \ f^+=\left(\ba{c}f_1\\f_2\\\vdots\\f_{n-1}\\f_n\ea\right),\ \ \
f^-=\left(\ba{c}f_n\\f_{n-1}\\\vdots\\f_2\\f_1\ea\right)\]
and define
\[U_n=\left(\Tni f^+,\,\dl^-\right)=\left(\Tni f^-,\,\dl^+\right).\]
If we set
\[\ph=1-U_n^2\]
then $\ph$ as a function of $t$ satisfies the equation
\bq\ph''=\hf\Big({1\ov\ph-1}+{1\ov\ph}\Big)(\ph')^2-{1\ov t}\ph'
-8\ph(\ph-1)+2{n^2\ov t^2}{\ph-1\ov\ph},\label{PV}\eq
which is a variant of the Painlev\'e V equation,\footnote{The
substitution $t^2=z$ and $\ph=w/(w-1)$ transforms (\ref{PV})
to the standard form of $P_V$ with parameters
$\alpha=0$, $\beta=-n^2/2$, $\gamma=2$ and $\delta=0$.}
 and in terms of
it $D_n$ has the
representation
\bq D_n(t)=\exp\left(4\int_0^t
\log(t/\tau)\,\tau\,\ph(\tau)\,d\tau\right).\label{Dform}\eq
This is reminiscent of the many representations now in the literature
for
Fredholm determinants in terms of Painlev\'e functions. We shall also
show that
\[W_n=U_n/U_{n-1}\]
satisfies
\bq W_n''={1\ov W_n}\,(W_n')^2-{1\ov t}\,W_n'+4\,{n-1\ov
t}\,W_n^2-{4n\ov
t}+4\,W_n^3-{4\ov W_n},
\label{PIII}\eq
which is a special case of the Painlev\'e III equation.

An important ingredient in the proofs is the following recurrence
relation satisfied
by the $U_n$:
\bq{n\ov t}\,U_n+(1-U_n^2)\,(U_{n-1}+U_{n+1})=0.\label{Urec}\eq
We shall see that this recurrence formula, 
sometimes known as the 
discrete Painlev\'e II equation (see, e.g.~\cite{Ni}), is equivalent
to one first shown to hold for  $\cU(n)$ by
V.~Periwal and D.~Shevitz~\cite{Pe}.  It was rediscovered
by M.~Hisakado \cite{H}, who also derived an equation
equivalent to (\ref{PV}) and
observed that this was one of the
Painlev\'e V equations which, by results of K.~Okamoto \cite{O}, is
reducible to Painlev\'e III.
Carrying through the Okamoto procedure is what led to our $W_n$,
although the proof we
give is direct. Our derivations of (\ref{PV}) and (\ref{Urec}) are
different from
those in \cite{H} and perhaps more
down-to-earth since we use only the simplest properties of Toeplitz
matrices and some
linear algebra. (They cannot be entirely unrelated, though, since the
orthogonal polynomials
which are central to the argument of \cite{H} can be defined in terms of
Toeplitz matrices.)

A remarkable connection between $\cU(n)$ and combinatorics was
discovered by
Gessel \cite{G}. Place the uniform measure on the symmetric group \S,
denote by
$\ell_N(\sigma)$ the length of the longest increasing subsequence in
$\sigma$, and
define $f_{Nn}$ by
\[{\rm Prob}\,(\ell_N(\sigma)\le n) ={f_{Nn}\ov N!}\, .\]
Then $D_n(t)$ is a generating function
for the $f_{Nn}$.\footnote{Gessel in \cite{G} does not write down the
symbol of
the Toeplitz matrix, nor does he
mention random matrices. But in light of the well-known
formula (\ref{ev}) and the subsequent work of Odlyzko {\it et al.}
\cite{Od} and
Rains \cite{R}, we believe it is fair to say that the connection
with random matrix theory begins with the discovery of (\ref{Dgen}).}
In fact
\bq D_n(t)=\sum_{N\ge 0} f_{Nn} {t^{2N}\over (N!)^2}\, . \label{Dgen}\eq
Recently, E.~Rains \cite{R} gave an elegant proof that
\bq f_{Nn}=\E\left(\vert\tr(U)\vert^{2N}\right),\label{probrep}\eq
which can be shown to be equivalent to (\ref{Dgen}) by a simple
argument.
Using the relationship (\ref{Dgen})  a sharp asymptotic result for the
distribution
function of the
random variable $\ell_N(\sigma)$ was recently obtained by J.~Baik,
P.~Deift
and K.~Johansson
\cite{BDJ}.
And at the same time they discovered yet another connection with
Painlev\'e.

Their main
result, which was quite difficult, was an asymptotic formula for
$D_n(t)$ which we now
describe. Introduce another parameter $s$ and suppose that $n$ and $t$
are related by
$n=2t+s\,t^{1/3}$.  Then as $t\ra\iy$ with fixed $s$ one has
\bq \lim_{t\ra\iy}e^{-t^2}\,D_{2t+s\,t^{1/3}}(t)=F(s).\label{Dasym}\eq
Here $F$ is the distribution function defined by
\bq F(s)=\exp\left(-\int_s^{\iy}(x-s)\,q(x)^2\,dx\right)\label{F}\eq
where $q$ is the solution of the Painlev\'e II equation
\bq q''=sq+2q^3\label{PII}\eq
satisfying $q(s)\sim \rm{Ai}(s)$ as $s\ra\iy$.  (For a proof
that such
a solution exists, see, e.g.~\cite{DZ}.)
Using a ``de-Poissonization'' lemma  due to Johansson \cite{J} these
asymptotics
for the generating function $D_n(t)$ led to the asymptotic formula
\bq\lim_{N\ra\iy}{\rm Prob}\,
\left({{\ell_N(\sigma)-2\sqrt{N}\over N^{1/6}}\leq s}\right)
=F(s)\label{dislim}\eq
for the distribution function of the normalized
 random variable $(\ell_N(\sigma)-2\sqrt{N})/N^{1/6}$.

It is a remarkable fact that this same
distribution function $F$ was first encountered
in random matrix theory where it gives the limiting
distribution for the normalized
largest eigenvalue in the Gaussian Unitary Ensemble
of Hermitian
matrices. More precisely, we have for this ensemble \cite{TW},
\[\lim_{N\ra\iy}{\rm Prob}\,\left(\left(\lambda_{\rm{max}}(N)-
\sqrt{2N}\,\right) \sqrt2\,N^{1/6}\leq s\right)=F(s).\]

In connection with these results just described, we shall do two things.
We show, first, how
one might have guessed the asymptotics (\ref{dislim}). More precisely,
we present a
simple argument
that if there is any limit theorem of this type, with $F(s)$ some
distribution
function and with some power $N^{\al}$ replacing $N^{1/6}$, then
necessarily $\al=1/6$ and $F$ is given by (\ref{F}) with $q$ a solution
to (\ref{PII}).
(The boundary condition on $q$, however, cannot be anticipated.) This
conclusion is
arrived at by considering the implications of (\ref{Dasym}) with
$t^{1/3}$ replaced
by $t^{2\al}$ for the recurrence formula (\ref{Urec}).

Secondly, we derive analogues of (\ref{probrep}) and (\ref{Dgen}) for
the subgroup $\cO_N$  of ``odd''
permutations of \S.\footnote{Our terminology for $\cO_N$ is
not standard.  For $N=2k$ one usually denotes  $\cO_N$
by ${\cal{B}}_k$, the hyperoctahedral group of order $2^k k!$ which
is the centralizer of the reversal permutation in \S.  Elements
of ${\cal{B}}_k$ are commonly called {\it signed permutations\/}.
Similar remarks hold for $N=2k+1$.}
 These are described as follows: if $N=2k$ think of
\S as acting on the integers from $-k$ to $k$ excluding $0$, and if
$N=2k+1$ think of
\S as acting on the integers from $-k$ to $k$ including $0$. In both
cases
$\si\in{\cal{S}}_N$ is called {\it odd} if
$\si(-n)=-\si(n)$ for all $n$. The number of elements in the subgroup
$\cO_N$ of odd
permutations
equals $2^k\,k!$ in both cases. Therefore if $b_{Nn}$ equals the number
of permutations
in $\cO_N$ having no increasing subsequence of length greater than $n$,
\bq{\rm Prob}(\ell_N(\si)\le n)={b_{Nn}\ov 2^k\,k!},\label{oddprob}\eq
where the uniform measure is placed on $\cO_N$. Rains \cite{R}
proved identities analogous to (\ref{probrep}) for these probabilities.
Using these we
are able to find representations for the two generating functions
\bq G_n(t)=\sum_{k\ge 0} b_{2k\, n} {t^{2k}\over (k!)^2},\ \ \
H_n(t)=\sum_{k\ge 0} b_{2k+1\, n} {t^{2k}\over (k!)^2},\label{GH}\eq
analogous to the representation (\ref{Dgen}). (See Theorem 1 below.)
The same determinants $D_n(t)$ arise as before but in the representation
for
$H_n(t)$,
whose derivation uses the machinery developed in earlier sections,
the quantities $U_n$ also appear. Once the representations are
established we can use
(\ref{Dasym}) and Johansson's lemma to deduce the
asymptotics of (\ref{oddprob}). We show that as $N\ra\iy$ we have for
fixed $s$
\bq{\rm Prob}\left({\ell_N(\si)- 2\sqrt{N}\ov 2^{2/3} N^{1/6}}\leq s
 \right)\ra F(s)^2,\label{oddlim}\eq
 where $F(s)$ is as in (\ref{dislim}).\sp

 In  Table~1 we give some statistics of the distribution functions
 $F$ and $F_\mathcal{O}:=F^2$.  In  Figure~1 we graph their densities.

\begin{table}
\begin{center}

\begin{tabular}{|l|cccc|}\hline
Distr & $\mu$ & $\sigma$ & $S$ & $K$ \\  \hline
$F$ & -1.77109 & 0.9018 & 0.224 & 0.093 \\
$F_\mathcal{O}$ & -1.26332 & 0.7789 & 0.329 & 0.225 \\ \hline
\end{tabular}
\vspace{2ex}
\caption{The mean ($\mu$),  standard deviation ($\sigma$),
skewness ($S$) and  kurtosis ($K$) of $F$ and $F_\mathcal{O}:=F^2$.}
\end{center}
\end{table}
\begin{figure}
\vspace{-.4in}
\begin{center}
\resizebox{7cm}{6cm}{\includegraphics{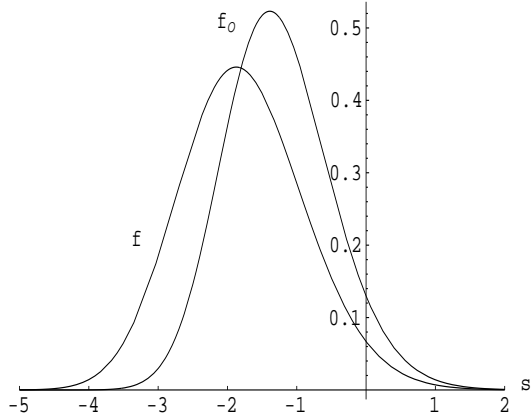}}
\vspace{2ex}
\caption{The probability densities $f=dF/ds$ and
$f_\mathcal{O}=dF_\mathcal{O}/ds$.}
\vspace{2ex}
\end{center}
\end{figure}

\setcounter{equation}{0}\renewcommand{\theequation}{2.\arabic{equation}}
\begin{center} {\bf II. The integral representation for {\boldmath
$D_n$}}\end{center}\sp

We write
\[\La=T_n(z\inv),\ \ \ \La'=T_n(z).\]
Thus $\La$ is the backward shift and $\La'$ is the forward shift. It is
easy to see that
\bq T_n(z\inv\,f)=T_n(f)\,\La+f^+\tn\dl^+=\La\,T_n(f)+\dl^-\tn
f^-,\label{Lcom}\eq
\bq T_n(z\,f)=T_n(f)\,\La'+f^-\tn\dl^-=\La'\,T_n(f)+\dl^+\tn
f^+,\label{L'com}\eq
where $\dl^{\pm}$ and $f^{\pm}$ were defined above
and $a\tn b$ denotes the matrix with $j,k$
entry $a_j\,b_k$. Relation (\ref{Lcom}) holds for any $f$ but
(\ref{L'com})
uses the fact that $f_{-k}=f_k$.

To derive (\ref{Dform}) we temporarily reintroduce variables $r$ and $s$
and set
$f(z)=e^{rz+sz\inv}$,
so $D_n$ and $T_n$ are functions of $r$ and $s$. Of course we areF
interested in $D_n(t,t)$.
We shall compute
\[\pl^2_{r,s}\,\log D_n(r,s)\Big|_{r=s=t}\]
in two different ways.

Using the fact that
\[\pl_s\,\log D_n(r,s)={\rm {tr}}\ T_n(f)\inv\, T_n(\pl_s f)=
{\rm {tr}}\ T_n(f)\inv\, T_n(z\inv f),\]
then differentiating with respect to $r$, we find that
\[\pl^2_{r,s}\,\log D_n(r,s)={\rm {tr}}\ [\Tni\,\Tn-
\Tni\, T_n(z\,f)\,\Tni \,T_n(z\inv\,f)]\]
\[={\rm {tr}}\ [I-\Tni\, T_n(z\,f)\,\Tni \,T_n(z\inv\,f)].\]
We now set $r=s=t$. Since $f$ is now as it was, we can use (\ref{Lcom})
and (\ref{L'com}).
If we multiply their first equalities on the left by $\Tni$ and use the
notation
$u^{\pm}=\Tni f^{\pm}$ we obtain
\[\Tni\, T_n(z\inv\,f)=\La+u^+\tn\dl^+,\ \ \ \Tni\,
T_n(z\,f)=\La'+u^-\tn\dl^-.\]
Hence the last trace equals that of
\[I-(\La'+u^-\tn\dl^-)(\La+u^+\tn\dl^+)\]\[=I-\La'\,\La-u^-\tn\La'\dl^-
-\La'u^+\tn\dl^+-(\dl^-,\,u^+)\,u^- \tn\dl^+.\]
The trace of $I-\La'\,\La$ equals 1, and
\[\La'\dl^-=0,\ \ \  {\rm {tr}}\
\La'u^+\tn\dl^+=(\La'u^+,\,\dl^+)=(u^+,\,\La\dl^+)=0,\]
so we have
\[\pl^2_{r,s}\,\log\,D_n(r,s)\Big|_{r=s=t}=1-(\dl^-,\,u^+)\,(u^-,\,\dl^+).\]
But $(u^-,\,\dl^+)=(\dl^-,\,u^+)=(\dl^-,\,\Tni f^+)=U_n$. Therefore
\bq \pl^2_{r,s}\,\log\,D_n(r,s)\Big|_{r=s=t}=1-U_n^2=\ph.\label{rs}\eq

Now let us go back to general $r$ and $s$. For any $p>0$ Cauchy's
theorem tells us that
the $j,k$ entry of $T_n(f)$ equals
\[{1\ov 2\pi i}\int_{|z|=p}e^{t\,(rz+sz\inv)}\,z^{-(j-k+1)}\,dz=p^{-j}\
{1\ov 2\pi i}\int_{|z|=1}e^{t\,(prz+p\inv sz\inv)}\,z^{-(j-k+1)}\,dz\
p^k.\]
It follows that $D_n(r,\,s)
=D_n(pr,\,s/p)$, and by analytic continuation this holds for any
(complex) $p$.
Setting $p=\sqrt{s/r}$ we see that
\[D_n(r,\,s)=D_n(\sqrt{rs},\,\sqrt{rs}).\]
It follows that $D_n(r,\,s)$ is a function of the product $rs$, as
stated
in the Introduction, and that
\[4\,\pl^2_{r,s}\,\log D_n(r,s)\Big|_{r=s=t}=
{d^2\ov dt^2}\log D_n(t,t)+{1\ov t}\,{d\ov dt}\log D_n(t,t).\]
Comparing this with (\ref{rs}) we see that we have shown
\[{d^2\ov dt^2}\log D_n(t,t)+{1\ov t}\,{d\ov dt}\log D_n(t,t)=4\,\ph.\]
This gives the representation (\ref{Dform}).

Of course it remains to show that this $\ph$ satisfies (\ref{PV}). We do
this by
first finding a formula for $dU_n/dt$ and then finding relations among
the various
quantities which occur for different values of $n$.\sp

\setcounter{equation}{0}\renewcommand{\theequation}{3.\arabic{equation}}
\begin{center} {\bf III. Differentiation}\end{center}\sp

In addition to $u^{\pm}=\Tni f^{\pm}$, we introduce $v^{\pm}=\Tni
\dl^{\pm}$ and we
compute some derivatives respect to $t$. First,
\[\dt\Tn=T_n((z+z\inv)f)=T_n(f)\,(\La+\La')+f^+\tn\dl^++f^-\tn\dl^-\]
by the first equalities of (\ref{Lcom}) and (\ref{L'com}), so
\bq\dt\Tni=-\Tni\,{d\Tn\ov dt}\,\Tni=-(\La+\La')\,\Tni-u^+\tn v^+-u^-\tn
v^-.\label{dTni}\eq
Next,
\[{df^+\ov \hspace{-1.5ex} dt}=
\left(\ba{c}f_0\\f_1\\\vdots\\f_{n-2}\\f_{n-1}\ea\right)+
\left(\ba{c}f_2\\f_3\\\vdots\\f_n\\f_{n+1}\ea\right)=\Tn\dl^++\La
f^++f_{n+1}\,\dl^-.\]
Hence
\bq\Tni\,{df^+\ov \hspace{-1.5ex}dt}=
\dl^++\Tni\La f^++f_{n+1}\,v^-.\label{Tf}\eq
Multiplying the second equality of (\ref{Lcom}) left and right by $\Tni$
gives
\bq \Tni\,\La=\La\,\Tni+u^+\tn v^+-v^-\tn u^-.\label{Tnicom}\eq
Therefore
\[\Tni\La f^+=\La\,u^++u^+\,(v^+,\,f^+)-v^-\,(u^-,\,f^+),\]
and substituting this into (\ref{Tf}) gives
\[\Tni\,{df^+\ov\hspace{-1.5ex}dt}=
\dl^++\La\,u^++u^+\,(v^+,\,f^+)+v^-\,\left(f_{n+1}-(u^-,\,f^+)\right).\]
Adding this to (\ref{dTni}) applied to $f^+$ gives
\[{du^+\ov\hspace{-1.5ex}dt}=
\dl^+-\La'\,u^++v^-\,(f_{n+1}-(u^-,\,f^+))-u^-\,(v^-,\,f^+).\]
Taking inner products with $\dl^-$ in the last displayed formula we
obtain (recall the
definition of $U_n$)
\[{dU_n\ov dt}=-(\La'\,u^+,\,\dl^-)+
(v^+,\,\dl^+)\,(f_{n+1}-(u^-,\,f^+))-(u^+,\,\dl^+)\,(v^+,\,f^-).\]
We used the fact, which follows from the symmetry of $\Tn$, that all our
inner products
whose entries have signs as superscripts are unchanged if both signs are
reversed.

To find $(\La'\,u^+,\,\dl^-)$, which is the same as
$(\La\,u^-,\,\dl^+)$, we observe that
\[\La\,f^-=\left(\ba{c}f_n\\f_{n-1}\\\vdots\\f_1\\0\ea\right)=
\Tn\,\dl^--f_0\,\dl^-.\]
Applying (\ref{Tnicom}) to $f^-$ therefore gives
\[\dl^--f_0\,v^-=\La\,u^-+u^+\,(v^+,\,f^-)-v^-\,(u^-,\,f^-).\]
Hence
\[-(\La'\,u^+,\,\dl^-)=-(\La\,u^-,\,\dl^+)=(u^+,\,\dl^+)\,(v^+,\,f^-)+
(v^-,\,\dl^+)\,(f_0-(u^+,\,f^+)).\]
Thus we have established the
differentiation formula
\bq{dU_n\ov dt}=
(v^-,\,\dl^+)\,\left(f_0-(u^+,\,f^+)\right)+
(v^+,\,\dl^+)\,\left(f_{n+1}-(u^-,\,f^+)\right).
\label{dum}\eq\sp

\setcounter{equation}{0}\renewcommand{\theequation}{4.\arabic{equation}}
\begin{center} {\bf IV. Relations}\end{center}\sp

New quantities appearing in the differentiation formula are
\[V_n^{\pm}=(v^{\pm},\,\dl^+).\]
There are others but we shall see that they may be expressed in terms of
these (with
different values of $n$), as indeed so will $U_n$.
We obtain our relations through several applications of the following
formula
for the inverse of a $2\times 2$ block matrix:
\bq\twotwo{A}{B}{C}{D}\inv=\twotwo{(A-BD\inv C)\inv}{\times}
{\times}{\times}.\label{ABCD}\eq
Here we assume $A$ and $D$ are square and the various inverses exist.
Only one block
of the inverse is displayed and the formula shows that $A-BD\inv C$
equals the inverse of
this block of the inverse matrix.
At first all that will be used about $f$ is
that $\Tn$ is symmetric. (There are modifications which hold in
general.)

We apply (\ref{ABCD}) first to the $(n+1)\times(n+1)$ matrix
\[\left(\begin{array}{cccc}0&0&\cd&1\\f_1&f_0&\cd&f_{n-1}\\
\vdots&\vdots&\cd&\vdots\\f_n&f_{n-1}&\cd&f_0\end{array}\right),\]
with $A=(0),\ \ D=\Tn,\ \ B=(0\;\cd\;0\;1),\ \ C=f^+$.
In this case $A-BD\inv C=-(\Tni\,f^+,\,\dl^-)=-U_n$. This equals the
reciprocal of the upper-left entry
of the inverse matrix, which in turn equals $(-1)^n$ times the
lower-left $n\times
n$ subdeterminant divided by $D_n$. Replacing the first
row by $(f_0\;f_1\;\cd\;f_n)$ gives the matrix
\[\left(\begin{array}{cccc}f_0&f_1&\cd&f_n\\f_1&f_0&\cd&f_{n-1}\\
\vdots&\vdots&\cd&\vdots\\f_n&f_{n-1}&\cd&f_0\end{array}\right)=T_{n+1}(f).\]
The upper-right
entry of its inverse equals on the one hand $V_{n+1}^-$ and on the other
hand $(-1)^n$ times the same subdeterminant as arose above divided by
$D_{n+1}$. This gives the identity
\bq -U_n=V_{n+1}^-\,{D_{n+1}\ov D_n}={V_{n+1}^-\ov
V_{n+1}^+}.\label{UV}\eq
(If we consider the polynomials on the circle which are
orthonormal with respect to the weight function $f$ then the right side
above is equal to
the constant term divided by the highest coefficient in the polynomial
of degree $n$.
Therefore our $-U_n$ equals the $S_{n-1}$ of \cite{H}.)

If we now take $A$ to be the upper-left corner of $T_{n+1}(f)$ and $D$
the complementary
$\Tn$ then $C=f^+$ and $B$ is its transpose, and we deduce that
\bq f_0-(u^+,\,f^+)={1\ov V_{n+1}^+}.\label{f0}\eq
To evaluate $f_{n+1}-(u^-,\,f^+)$, the other odd coefficient appearing
in (\ref{dum}),
we consider
\[\left(\begin{array}{ccccc}f_0&f_1&\cd&f_n&f_{n+1}\\
f_1&f_0&\cd&f_{n-1}&f_n\\
\vdots&\vdots&\cd&\vdots&\vdots\\
f_n&f_{n-1}&\cd&f_0&f_1\\
f_{n+1}&f_n&\cd&f_1&f_0\end{array}\right)=T_{n+2}(f).\]
We apply to this an obvious modification of (\ref{ABCD}), where $A$ is
the $2\times 2$ matrix consisting of the four corners of the large
matrix, $D$ is the
central $\Tn$, $C$ consists of the two columns $f^+$ and $f^-$ and $B$
consists of the
rows which are their transposes. Then
\[A-BD\inv
C=\twotwo{f_0-(u^+,\,f^+)}{f_{n+1}-(u^-,\,f^+)}{f_{n+1}-(u^-,\,f^+)}
{f_0-(u^+,\,f^+)}\]
and our formula tells us that this is the inverse of
\[\twotwo{V_{n+2}^+}{V_{n+2}^-}{V_{n+2}^-}{V_{n+2}^+}.\]
This gives the two formulas
\[f_0-(u^+,\,f^+)={V_{n+2}^+\ov {V_{n+2}^+}^2-{V_{n+2}^-}^2}\,,\ \ \
f_{n+1}-(u^-,\,f^+)={-V_{n+2}^-\ov {V_{n+2}^+}^2-{V_{n+2}^-}^2}\,.\]
Comparing the first with (\ref{f0}) we see that
\bq{V_{n+2}^+}^2-{V_{n+2}^-}^2=V_{n+1}^+\,V_{n+2}^+,\label{V}\eq
and therefore that the preceding relations can be written
\bq f_0-(u^+,\,f^+)={1\ov V_{n+1}^+},\ \ \
f_{n+1}-(u^-,\,f^+)=-{1\ov V_{n+1}^+}\,{V_{n+2}^-\ov
V_{n+2}^+}.\label{fn}\eq
Notice that (\ref{UV}) and (\ref{V}) give
\bq 1-U_n^2={V_n^+\ov V_{n+1}^+}.\label{UV1}\eq
This is our $\ph$.\sp

The relations we obtained so far in this section are completely general.
The recurrence
(\ref{Urec}), however, depends on our specific
function $f$. Integration by parts gives
\bq k\,f_k={t\ov2\pi i}\int(z-z\inv)\,e^{t\,(z+z\inv)}\,z^{-k-1}\,dz=
t\,(f_{k-1}-f_{k+1}).\label{fdiff}\eq
Hence if $M={\rm {diag}}\;(1\;2\;\cd\;n)$ we have
\[M\,\Tn-\Tn\,M=t\,T_n((z-z\inv)\,f),\]
and by the first identities of (\ref{Lcom}) and (\ref{L'com}) this
equals
\[t\,\Big[\Tn\,(\La'-\La)+f^-\tn\dl^--f^+\tn\dl^+\Big],\]
so
\[\Tni\,M-M\,\Tni=t\,\Big[(\La'-\La)\,\Tni+u^-\tn v^--u^+\tn v^+\Big].\]
Applying this to $\dl^-$ gives
\bq
n\,v^--M\,v^-=t\,\Big[(\La'-\La)\,v^-+u^-\,
(v^-,\,\dl^-)-u^+\,(v^+,\,\dl^-)\Big].
\label{vid}\eq
Now (\ref{fdiff}) say
\[M\,f^+=t\,\left(\ba{c}f_0-f_2\\f_1-f_3\\\vdots\\f_{n-2}-
f_n\\f_{n-1}-f_{n+1}\ea\right)\]
whereas (this is relevant since the transpose of $\La'-\La$ is
$\La-\La'$)
\[(\La-\La')\,f^+=\left(\ba{c}f_2\\f_3-f_1\\\vdots\\f_n-
f_{n-2}\\-f_{n-1}\ea\right).\]
Therefore
\[M\,f^++t\,(\La-\La')\,f^+=t\,(f_0\,\dl^+-f_{n+1}\,\dl^-),\]
and so taking inner products with $f^+$ in (\ref{vid}) gives
\[{n\ov t}\,(v^-,\,f^+)=
f_0\,(v^-,\,\dl^+)-f_{n+1}\,(v^-,\,\dl^-)+(u^-,\,f^+)\,(v^-,\,\dl^-)-
(u^+,\,f^+)\,(v^+,\,\dl^-)\]
or equivalently, since $(v^-,\,f^+)=U_n$,
\[{n\ov t}\,U_n=
\left(f_0-(u^+,\,f^+)\right)\,V_n^--\left(f_{n+1}-(u^-,\,f^+)\right)\,V_n^+.\]
Using (\ref{fn}) we rewrite this as
\[{n\ov t}\,U_n={V_n^-\ov V_{n+1}^+}+{V_n^+\ov V_{n+1}^+}\,{V_{n+2}^-\ov
V_{n+2}^+}
={V_n^+\ov V_{n+1}^+}\,{V_n^-\ov V_n^+}+{V_n^+\ov
V_{n+1}^+}\,{V_{n+2}^-\ov V_{n+2}^+}.\]
Using this, (\ref{UV1}) and (\ref{UV}) we arrive at (\ref{Urec}).\sp

\setcounter{equation}{0}\renewcommand{\theequation}{5.\arabic{equation}}
\begin{center} {\bf V. Painlev\'e V and Painlev\'e III}\end{center}\sp

We first show that $\ph$ satisfies (\ref{PV}). Our formula (\ref{dum})
for $dU_n/dt$
can now be written
\bq{dU_n\ov dt}={V_n^-\ov V_{n+1}^+}-{V_n^+\ov V_{n+1}^+}\,{V_{n+2}^-\ov
V_{n+2}^+}
={V_n^+\ov V_{n+1}^+}\,\Big({V_n^-\ov V_n^+}-{V_{n+2}^-\ov
V_{n+2}^+}\Big)
=-(1-U_n^2)\,(U_{n-1}-U_{n+1}),\label{dU1}\eq
by (\ref{fn}), (\ref{UV1}) and (\ref{UV}). Adding and subtracting
(\ref{Urec}) gives
us the two formulas
\bq {dU_n\ov dt}={n\ov t}\,U_n+2\,U_{n+1}\,(1-U_n^2),\label{Uder1}\eq
\bq{dU_n\ov dt}=-{n\ov t}\,U_n-2\,U_{n-1}\,(1-U_n^2).\label{Uder2}\eq
These are equations (4.5) and (4.6) of \cite{H}. As was done there,
we solve (\ref{Uder1}) for $U_{n+1}$ in terms of $U_n$ and $dU_n/dt$
and substitute this into (\ref{Uder2}) with $n$ replaced by $n+1$. We
get a
second-order differential equation for $U_n$ which is equivalent to
equation (\ref{PV}) for $\ph=1-U_n^2$.\sp

Next we show that $W_n=U_n/U_{n-1}$ satisfies (\ref{PIII}).
In computing the derivative of $W_n$ we use (\ref{Uder2})
to compute the derivative of $U_n$ and (\ref{Uder1}) with $n$ replaced
by $n-1$ to compute
the derivative of $U_{n-1}$. We get
\bq W_n'=-{2n-1\ov t}\,W_n-2+4\,U_n^2-2\,W_n^2.\label{W}\eq
Using (\ref{Uder2}) once again we compute
\[(U_n^2)'=2\,U_n\,\left(-{n\ov t}\,U_n-2\,U_{n-1}\,(1-U_n^2)\right)=
-2\,{n\ov t}\,U_n^2-4\,{U_n^2\,(1-U_n^2)\ov W_n}.\]
Differentiating (\ref{W}) and using this expression for $(U_n^2)'$ we
obtain a
formula for
$W_n''$ in terms of $W_n,\ W_n'$ and $U_n^2$. Then we solve (\ref{W})
for $U_n^2$
in terms of $W_n$ and $W_n'$. Substituting this into the formula for
$W_n''$ gives
(\ref{PIII}).\sp

In order to specify the solutions of the equations (\ref{PV}) and
(\ref{PIII}) we must
determine the initial conditions at $t=0$. Clearly $\ph(0)=1$ but this
does not determine
$\ph$ uniquely. One can see that $\ph^{(k)}(0)=0$ for $k<2n$ and that
what determines $\ph$
uniquely is $\ph^{(2n)}(0)$. We shall show that
\bq\ph^{(2n)}(0)=-{(2n)!\ov n!^2}.\label{init}\eq

By (\ref{UV})
$U_n=V_{n+1}^-/V_{n+1}^+$. Now $V_{n+1}^+$ is the upper-left corner of
$T_{n+1}(f)^{-1}$ and so
tends
to 1 as $t\rightarrow 0$. So let us see how $V_{n+1}^-$, which is the
upper-right corner
of $T_{n+1}(f)^{-1}$, behaves. More exactly, let us find the term in its
expansion with the lowest power
of $t$.

We have
\[e^{2t\cos\th}=\sum_k{t^k\ov k!}\,(e^{i\th}+e^{-i\th})^k\]
\[=\sum_{0\leq j\leq k}{t^k\ov k!}\,C(k,j)\,e^{-i(2j-k)\th}
=\sum_{|j|\leq k}{t^k\ov k!}\,C(k,(j+k)/2)\,e^{-ij\th}.\]
This gives
\bq \Tn=\sum_{|j|\leq k}{t^k\ov k!}\,C(k,(j+k)/2)\,\La^j
=I+\sum_{k>0 \atop |j|\leq k}{t^k\ov
k!}\,C(k,(j+k)/2)\,\La^j.\label{sum}\eq
(Here $\La^j$ denotes the usual power when $j\geq0$, but when $j<0$ it
denotes ${\La'}^{|j|}$.)

We use the Neumann expansion
\[\Tni=I+\sum_{l\geq1}(-1)^l\left(\sum_{k>0 \atop |j|\leq k}{t^k\ov
k!}\,C(k,(j+k)/2)\,\La^j\right)^l.\]
If we expand this out we get a sum of terms of the form coefficient
times
\[t^{k_1+\cdots+k_l}\,\La^{j_1}\,\cdots\,\La^{j_l}.\]
Now the product $\La^{j_1}\,\cdots\,\La^{j_l}$ can
only have a nonzero upper-right entry when $j_1+\cdots+j_l\geq n$. Since
each
$|j_i|\leq k$ the power of $t$ must be at least $n$, and this power
occurs only when
each $j_i=k$. That means that we get the same lowest power of $t$ term
for
the upper-right entry if in (\ref{sum})
we only take the terms with $j=k$, in other words of we replace
$T_{n+1}(f)$ by
\[\sum_{k\geq 0}{t^k\ov k!}\,\La^k=e^{t\,\La}.\]
The inverse of this operator is $e^{-t\,\La}$ and the upper-right corner
of this matrix
is exactly $(-1)^n\,t^n/n!$. This shows that
\[ U_n=-V_{n+1}^-/V_{n+1}^+=(-1)^{n+1}\,t^n/n!+O(t^{n+1}),\]
and so
\[\ph=1-U_n^2=1-{t^{2n}\ov(n!)^2}+O(t^{2n+1}),\]
which gives (\ref{init}). We also see that $W_n=U_n/U_{n-1}$ satisfies
the initial condition
\[W_n(t)=-{t\ov n}+O(t^2).\]
Using the differential equation (\ref{PIII}) together with this initial
condition we find
\[ W_n(t)=-{t\ov n}-{t^3\ov n^2(n+1)}-{2t^5\over n^3(n+1)(n+2)}-
{(5n+6)t^7\ov n^4(n+1)^2(n+2)(n+3)}+O(t^9).\]\sp

\setcounter{equation}{0}\renewcommand{\theequation}{6.\arabic{equation}}
\begin{center} {\bf VI. Painlev\'e II}\end{center}\sp

We present a heuristic argument that if there is any limit theorem of
the type
(\ref{dislim}), with some distribution
function $F(s)$ and some power $N^{\al}$ replacing $N^{1/6}$, then
necessarily $\al=1/6$ and $F$ is given by (\ref{F}). First we note that
Johansson's lemma (which we shall state in the next section) leads from
(\ref{dislim}) to
(\ref{Dasym}) with the power
$t^{1/3}$ replaced by $t^{2\al}$. We assume that
$F$ is smooth and that the limit in (\ref{Dasym})
commutes with $d/ds$, so that taking the second logarithmic derivative
gives
\[\lim_{t\ra\iy}{d^2\ov ds^2}\log\,D_{2t+s\,t^{2\al}}=-q(s)^2,\]
where $q^2$ is now {\it defined}  by $-q^2=(\log\,F)''$ and $q$ is
defined to be
the positive square root of $q^2$ (for large $s$).

Since changing $n=2t+s\,t^{2\al}$ by 1 is the same as changing $s$ by
$t^{-2\al}$,
we have the large $t$ asymptotics
$$\log\,D_{n+1}+\log\,D_{n-1}-2\,\log\,D_{n}\sim
t^{-4\al}\,{d^2\ov dt^2}\log\,D_{2t+s\,t^{2\al}}\sim
-t^{-4\al}\,q(s)^2.$$
On the other hand $V_n^+=D_{n-1}/D_n$ and so
\bq {D_{n+1}\,D_{n-1}\ov D_{n}^2}={V_n^+\ov
V_{n+1}^+}=1-U_n^2,\label{UVD}\eq
by (\ref{UV1}). We deduce
\bq \log\,(1-U_n^2)\sim -t^{-4\al}\,q(s)^2,\ \ \ U_n^2\sim
t^{-4\al}\,q(s)^2,\label{Uasym} \eq

Now the $U_n$ are of variable sign, as is clear from (\ref{Urec}). Let
us consider
those $n$ going to infinity such that
\bq U_{n-1}\ge0,\ \ U_n\le0,\ \ U_{n+1}\ge0,\label{Usign}\eq
and write (\ref{Urec}) as
\bq
t\,(U_{n+1}+U_{n-1}+2\,U_n)\,(1-U_n^2)=-(n-2t)\,U_n-2t\,U_n^3.\label{UPII}\eq
Because of (\ref{Usign}), (\ref{Uasym}) and the fact that changing
$n=2t+s\,t^{2\al}$ by
1 is the same as changing $s$ by $t^{-2\al}$, we have when $t$ is large,
\[U_{n+1}+U_{n-1}+2\,U_n \sim t^{-6\al}\,q''(s).\]
Since also
\[n-2t \sim s\,t^{2\al},\ \ \ U_n\sim -t^{-2\al}\,q(s),\]
(\ref{UPII}) becomes the approximation
\bq t^{1-6\al}\,q''(s)\approx
s\,q(s)+2\,t^{1-6\al}\,q(s)^3.\label{qeq}\eq

Let us show that $\al=1/6$. If $\al>1/6$ then letting
$t\ra\iy$ in (\ref{qeq}) gives $q(s)=0$ and so $F$ is the exponential of
a linear function
and therefore not a distribution function.
If $\al<1/6$ then dividing by $t^{1-6\al}$ and letting $t\ra\iy$ in
(\ref{qeq})
gives $q''(s)=2\,q^3$. Solving this gives two sets of solutions
\[s=\pm\int_q^{q_0}{dq\ov\sqrt{q^4+c_1}}+c_2,\]
where $q_0,\ c_1$ and $c_2$ can be arbitrary. Now $q$ is small when $s$
is
large and positive, so $s$ is large and positive when $q$ is small.
Therefore we have
to have the + sign and we must have $c_1=0$.
Then $F(s)$ is of the form $|s-c|^{-1}$ times the exponential of a
linear function and
therefore is not a distribution function. The only remaining case is
$\al=1/6$, and
then (\ref{qeq}) becomes (\ref{PII}). It follows that $F(s)$ must be
given by (\ref{F})
times the exponential of a linear function. This extra factor must be 1
since
(\ref{F}) is already a distribution function.

Now to derive this we assumed that the $n$ under consideration were such
that
(\ref{Usign}) held. We would have reached the same conclusion if all the
inequalities were
reversed. If $n\ra\iy$ in such a way that, say, $U_{n-1}$ and $U_n$ have
one sign and
$U_{n+1}$ the other, then (the reader can check this) we would have
reached the
conclusion $q=0$. Thus the only possibility for $q$ to give a
distribution
function occurs when $\al=1/6$ and $q$ satisfies (\ref{PII}).\sp

\setcounter{equation}{0}\renewcommand{\theequation}{7.\arabic{equation}}
\begin{center} {\bf VII. Odd permutations}\end{center}\sp

Recall that $b_{Nn}$ equals the number of permutations in $\cO_N$
having no increasing subsequence of length greater than $n$. The
representations of Rains
\cite{R} for these quantities are
\bq b_{2k\, n}=\E\left(\vert\tr\,(U^2)^k\vert^2\right),\label{R1}\eq
\bq b_{2k+1\,
n}=\E\left(\vert\tr\,(U^2)^k\,\tr\,(U)\vert^2\right).\label{R2}\eq
\sp

\noi{\bf Theorem 1}. Let $G_n(t)$ and $H_n(t)$ be the generating
functions defined in (\ref{GH}).
Then
\begin{eqnarray}
\label{Gn}G_n(t)&=&\left\{
        \begin{array}{ll}
        D_{{n\over 2}}(t)^2, &\hspace{1.2em}  n\ {\rm even,} \\&\\
        D_{{n-1\over 2}}(t) D_{{n+1\over 2}}(t), & \hspace{1.2em} n\
{\rm odd,}
        \end{array}\right.\\
&&\nonumber\\
\label{Hn}H_n(t)&=&\left\{
        \begin{array}{ll}
        D_{{n\over 2}-1}(t) D_{{n\over 2}+1}(t),
        &\hspace{1em}n\ {\rm even,} \\&\\
 D_{{n-1\over 2}}(t) D_{{n+1\over 2}}(t), & \hspace{1em} n\ {\rm odd.}
        \end{array}\right.
\end{eqnarray}
\sp

We prove a lemma which gives a preliminary representation for the
generating function
in terms of other Toeplitz determinants. Let
\[ g(z,t_1,t_2)=
g(z)=e^{t_1(z+z\inv)+t_2(z^2+z^{-2})}\]
and define
\[ \hat{D}_n(t_1,t_2)=\hat{D}_n=\det T_n(g). \]\sp

\noi{\bf Lemma}. We have
\bq G_n(t)=\hat{D}_n(0,t)\label{G2}\eq
\bq H_n(t)={1\over 4} {\partial^2 \hat{D}_n\over \partial t_1^2}(0,t)
+{1\over 4} {\partial^2 \hat{D}_n\over \partial
t_1^2}(0,-t).\label{H2}\eq\sp

The proof of (\ref{G2}) is essentially the same as the proof that
(\ref{Dgen}) and
(\ref{probrep}) are equivalent. First observe that
\[ \E\left(\left(\tr(U^2)+\overline{\tr(U^2)}\right)^{2k}\right)
= \sum_{m=0}^{2k}{2k\choose m}
\E\left(\tr(U^2)^m\overline{\tr(U^2)}^{2k-m}\right).\]
Each summand with $m\neq k$ vanishes since by the invariance of Haar
measure replacing
each $U$ by $\zeta U$, with $\zeta$ a complex number of absolute value
1, does not
change the
summand but at the same time multiplies it by $\zeta^{4m-4k}$. Thus,
\[\E\left(\left(\tr(U^2)+\overline{\tr(U^2)}\right)^{2k}\right)
= {2k\choose k} \E\left(\tr(U^2)^k\overline{\tr(U^2)}^k\right).\]
Hence (\ref{R1}) is equivalent to
\[b_{2k\, n}={(k!)^2\ov (2k)!}\,
\E\left(\left(\tr(U^2)+\overline{\tr(U^2)}\right)^{2k}\right).\]
Therefore if the eigenvalues of $U$ are $e^{i\th_1},\cd,e^{i\th_n}$ we
have
\[G_n(t)=\sum_{k\ge 0} b_{2k\, n} {t^{2k}\over (k!)^2}
=\sum_{k\ge 0} \E\left(\left(\sum\cos
2\theta_j\right)^{2k}\right){(2t)^{2k}\over (2k)!}\]
\[=\E\left(\prod_{j=1}^n e^{2t\cos 2\theta_j}\right)=\hat{D}_n(0,t).\]
The last step follows from (\ref{ev}).

This gives (\ref{G2}). To prove (\ref{H2}) we use (\ref{R2}):
\[b_{2k+1\,
n}=\E\left(\tr(U)\tr(U^2)^k\overline{\tr(U)}\overline{\tr(U^2)}^k\right)\]
\[={1\ov 2}{(k!)^2\ov (2k)!}\,
\E\left(\left(\tr(U)+\overline{\tr(U)}\right)^2
\left(\tr(U^2)+\overline{\tr(U^2)}\right)^{2k}\right),\]
by expanding the right side as before. Hence
\[H_n(t)={1\over 2}\sum_{k\ge 0} {t^{2k}\over (2k)!}
\E\left(\left(\tr(U)+\overline{\tr(U)}\right)^2
\left(\tr(U^2)+\overline{\tr(U^2)}\right)^{2k}\right)\]
\[=2\, \E\left(\left(\sum\cos\theta_j\right)^2
\cosh\left(2t\sum\cos\theta_j\right)\right)\]
\[=\E\left(\left(\sum\cos\theta_j\right)^2\prod_{j=1}^n
 e^{2t\cos 2\theta_j}
\right) +\E\left(\left(\sum\cos\theta_j\right)^2\prod_{j=1}^n
 e^{-2t\cos 2\theta_j}\right)\]
 \[={1\over 4} {\partial^2 \over \partial t_1^2}\E\left(\prod_{j=1}^n
 g(e^{i\theta_j},t_1,t_2)\right)(0,t)
 +{1\over 4} {\partial^2 \over \partial t_1^2}\E\left(\prod_{j=1}^n
 g(e^{i\theta_j},t_1,t_2)\right)(0,-t).\]
 From the last equality (\ref{H2}) follows.\sp

 To prove the theorem we consider (\ref{Gn}) first. Observe that $\hat
D_n(0,t)$ is
 the determinant of $T_n(h)$
 where $h(z)=f(z^2)$. It has Fourier coefficients
 \[ h_{2k}=f_k,\ \ \ h_{2k+1}=0.\]
 Let us rearrange the basis vectors $e_0,\;e_1,\cd,e_{n-1}$ of
 our underlying $n$-dimensional space as
 \bq e_0,\;e_2,\cd,\ \ e_1,\;e_3,\cd.\label{vecs}\eq
 Then we see from the above that $T_n(h)$ becomes the direct sum of two
Toeplitz matrices
 associated with $f$, the orders of these matrices being
 the sizes of the two groups of basis vectors in (\ref{vecs}). If $n$ is
even both
 groups have size $n/2$ whereas if $n$ is odd the sizes are \linebreak
$(n\pm 1)/2$.
 Since $\hat D_n(0,t)=\det T_n(h)$ is the product of the
 corresponding Toeplitz determinants associated with $f$, we have
(\ref{Gn}).

 The proof of (\ref{Hn}) is not so simple. We have
 \[{1\over\hat{D}_n}{\partial^2\hat{D}_n\over\pl t_1^2}=
 \pl_{t_1}^2\log\hat{D}_n+\left(\pl_{t_1}\log\hat{D}_n\right)^2.\]
 Now, as in the computation leading to (\ref{rs}),
\[ \pl_{t_1}^2\log\hat{D}_n(t_1,t_2)
=\tr\left(T_n(g)^{-1}\,T_n(\pl_{t_1}^2 g)-T_n(g)^{-1}\,T_n(\pl_{t_1}g)
\,T_n(g)^{-1}
\,T_n(\pl_{t_1}g)\right)\]
\[=\tr\left[T_n(g)^{-1}\,T_n((z+z\inv)^2 g)-T_n(g)^{-1}\,T_n((z+z\inv)g)
\,T_n(g)^{-1}
\,T_n((z+z\inv)g)\right].\]
This is to be evaluated first at $t_1=0,\ t_2=t$.
Since $g(z,0,t)=f(z^2)=h(z)$ we must compute
\[\tr\left[T_n(h)^{-1}\,T_n((z+z\inv)^2 h)-T_n(h)^{-1}\,T_n((z+z\inv)h)
\,T_n(h)^{-1}
\,T_n((z+z\inv)h)\right].\]

We write $w^{\pm}=T_n(h)\inv h^{\pm}$, so that the $w^{\pm}$ are
associated with $h$
just as $u^{\pm}$ are associated with $f$. Using (\ref{Lcom}) and
(\ref{L'com}) we find
that
\[T_n(h)^{-1}\,T_n((z+z\inv)h)\, T_n(h)^{-1}\,T_n((z+z\inv)h)=
(\La'+w^-\otimes \dl^-+\La+w^+\otimes\dl^+)^2\]
and from this that
\[\tr\left(T_n(h)^{-1}\,T_n((z+z\inv)h)
\,T_n(h)^{-1}\,T_n((z+z\inv)h)\right)=\]
\[2n-2+2(w^-,\La\dl^-)+2(w^+,\La'\dl^+)+2(\dl^-,w^+)\,(w^-,\dl^+)+
(\dl^-,w^-)^2+(\dl^+,w^+)^2.\]
Consequently
\[\tr\left[T_n(h)^{-1}\,T_n((z+z\inv)^2 h)-T_n(h)^{-1}\,T_n((z+z\inv)h)
\,T_n(h)^{-1}
\,T_n((z+z\inv)h)\right]\]
\[=\tr\left(T_n(h)^{-1}\,T_n((z^2+z^{-2})h)\right)\]
\[+2-2(w^-,\La\dl^-)-2(w^+,\La'\dl^+)-2(\dl^-,w^+)\,(w^-,\dl^+)-
(\dl^-,w^-)^2-(\dl^+,w^+)^2.\]
This is $\pl_{t_1}^2\log\hat{D}_n(0,t)$. Similarly, we find
\[\pl_{t_1}\log\hat{D}_n(0,t)=
 (w^+,\delta^+)+(w^-,\delta^-),\]
 so that when $t_1=0,\ t_2=t$
 \[{1\ov\hat{D}_n}{\pl^2\hat{D}_n\ov\pl t_1^2}
 =\pl_{t_1}^2\log\hat{D}_n+\left(\pl_{t_1}\log\hat{D}_n\right)^2
=\tr\left(T_n(h)^{-1}\,T_n((z^2+z^{-2})h)\right)\]
\[+2-2(w^-,\La\dl^-)-2(w^+,\La'\dl^+)
 -2(w^+,\dl^-)\,(w^-,\dl^+)+2(w^+,\dl^+)\,(w^-,\dl^-).\]
 Since $T_n(h)$ is symmetric all superscripts in the symbols in the
inner product may
 be reversed as long as we interchange $\La$ and $\La'$. We
 therefore  have shown
 \[{1\ov\hat{D}_n(0,t)}{\pl^2\hat{D}_n(0,t)\ov\pl t_1^2}=
 \tr\left(T_n(h)^{-1}\,T_n((z^2+z^{-2})h)\right)
 +2-4\,(w^+,\La'\dl^+)-2\,(w^+,\dl^-)^2+2\,(w^+,\dl^+)^2.\]

 Let us rearrange our basis elements as in (\ref{vecs}) and suppose the
first group
 has $n_1$ vectors and the second group has $n_2$. Then $T_n(h)\inv$
becomes the
 matrix direct sum
 \[\twotwo{T_{n_1}(f)\inv}{0}{0}{T_{n_2}(f)\inv}\]
 and $T_n(h)^{-1}T_n\,((z^2+z^{-2})h)$ becomes
 \[\twotwo{T_{n_1}(f)\inv \,T_{n_1}((z+z\inv)f)}{0}{0}{T_{n_2}(f)\inv\,
T_{n_2}((z+z\inv)f)}.\]
By a now familiar computation we find from this that
\[\tr\left(T_n(h)^{-1}\,T_n((z^2+z^{-2})h)\right)=2\,(u^+_{n_1},\,\dl^+_{n_1})
+2\,(u^+_{n_2},\,\dl^+_{n_2}),\]
where $u^+_m$ and $\dl^+_m$ denote the quantities $u^+$ and $\dl^+$
associated with the
index $m$. We use similar notation below. To continue, after rearranging
our basis we
have the replacements
\[h^+\ra\twoone{0}{f^+_{n_2}},\ \ \ w^+\ra\twoone{0}{u^+_{n_2}},\ \ \
\dl^+\ra\twoone{\dl^+_{n_1}}{0},\ \ \
\La'\dl^+\ra\twoone{0}{\dl^+_{n_2}}\]
and
\[\dl^-\ra\twoone{0}{\dl^-_{n_2}}\ {\rm if}\ n\ {\rm is\ even},\ \ \ \ \
\dl^-\ra\twoone{\dl^-_{n_1}}{0}\ {\rm if}\ n\ {\rm is\ odd}.\]
  It follows from these that
 \[(w^+,\,\dl^+)=0,\ \ (w^+,\,\La'\dl^+)=(u^+_{n_2},\,\dl^+_{n_2})\]
 and that $(w^+,\,\dl^-)=(u^+_{n_2},\,\dl^-_{n_2})$ if $n$ is even and
 $(\dl^-,\,w^+)=0$ if $n$ is odd.

If we modify our notation by writing
 \[U_n^{\pm}=(u_n^+,\,\dl_n^{\pm}),\]
so that $U_n^-$ is what we have been denoting by $U_n$,
 the above gives
 \[{1\ov 2}{1\ov\hat{D}_n(0,t)}\,{\pl^2\hat{D}_n(0,t)\ov\pl t_1^2}
 =1+U^+_{n_1}-U^+_{n_2}-(w^+,\,\dl^-)^2=
 \left\{\begin{array}{ll}1-{U^-_{{n\ov2}}}(t)^2,&n\ {\rm even},\\&\\
 1+U^+_{{n+1\ov2}}(t)-U^+_{{n-1\ov2}}(t), & n\ {\rm
odd}.\end{array}\right.\]

 To evaluate our quantities at $(0,-t)$ we observe that if $C$ is the
diagonal matrix
 with diagonal entries $1,\,-1,\, 1,\cd,(-1)^n$ and we replace $t$ by
$-t$ then we have the
 replacements
 $T_n(f)\ra CT_n(f)C$ and $f^+\ra-C f^+$ and therefore $u^+\ra-C u^+$.
Therefore
 also $U_n^+=(u^+,\,\dl^+)\ra-U_n^+$ and
 $U_n^-=(u^+,\,\dl^-)\ra(-1)^{n+1}\,U_n^-$. Hence in the last
 displayed formula ${U^-_{n/2}}^2$ is an
 even function of $t$ whereas $U^+_{(n\pm1)/2}$ are odd functions
 of $t$. Also, $\hat{D}_n(0,t)$ is an even function of $t$. Thus
 \[H_n(t)=G_n(t)\times \left\{
 \begin{array}{ll}
       1-U^-_{{n\ov 2}}(t)^2,
        &\> n\ {\rm even,} \\&\\
 1, & \> n\ {\rm odd.}
        \end{array}\right.\]
Recalling (\ref{UV1}) and the general fact $V^+_m(t)=D_{m-1}(t)/D_m(t)$
we
obtain (\ref{Hn}). \sp

 \noi{\bf Theorem 2}. Let $\ell_N(\sigma)$ denote the length
 of the longest increasing subsequence of $\sigma$ in the subgroup
$\cO_N$ of
 \S. Then
\[\lim_{N\ra\iy}{\rm Prob}\left({\ell_N(\si)- 2\sqrt{N}\ov 2^{2/3}
N^{1/6}}\leq s\right)=F(s)^2,\]
 where $F(s)$ is as in (\ref{dislim}).\sp

For the proof we shall apply Johansson's lemma \cite{J}, which we now
state:\sp

\noi {\bf Lemma}. Let $\{P_k(n)\}_{k\ge 0}$ be
a family of distribution functions defined on the nonnegative integers
$n$
and $\varphi_n(\lambda)$ the generating function
\[ \vp_n(\lambda)= e^{-\lambda} \sum_{k\ge 0} P_k(n) {\lambda^k\over
k!}.\]
(Set $P_0(n)=1$.)
Suppose that for all $n$, $k\ge 1$,
\bq P_{k+1}(n)\le P_k(n). \label{monotone} \eq
If we define $\mu_k=k+4\sqrt{k\log k}$ and $\nu_k=k-4\sqrt{k\log k}$,
then
there is a constant $C$ such that
\[\vp_n(\mu_k) - {C\over k^2} \le P_k(n) \le \vp_n(\nu_k) +{C\over
k^2}\]
for all sufficiently large $k$, $0\le n \le k$.\sp

This allows one to deduce that $P_k(n)\sim \vp_n(k)$ as $n,\, k\ra\iy$
under
suitable conditions,
and is how one obtains the equivalence of (\ref{Dasym}) and
(\ref{dislim}). We shall
apply the lemma to the distributions functions
\[\varphi^e_n(\la)=e^{-\la} G_n\left(\sqrt{\la/2}\right)
=e^{-\la} \sum_{k\ge 0} F_{2k}(n) {\lambda^k\over k!}\, ,\]
\[\varphi^o_n(\la)=e^{-\la} H_n\left(\sqrt{\la/2}\right)
=e^{-\la} \sum_{k\ge 0} F_{2k+1}(n) {\lambda^k\over k!}\]
where
\bq F_N(n)={\rm Prob}(\ell_N(\si)\le n)={b_{Nn}\ov 2^k\,k!}\label{FN}\eq
when $N=2k$ or $2k+1$.
To apply this lemma we must prove another\sp

\noi{\bf Lemma}. We have $F_{N+2}\leq F_N$ for all $N$.\sp

We show this simultaneously for $N=2k$ and $N=2k+1$. Take a
$\si\in\cO_{N+2}$ and
remove the two-point set $\{-1,\,1\}$ from its domain. Then $\si$ maps
the remaining
$N$-point set one-one
onto another $N$-point set. If we identify both of these sets with the
integers from $-k$
to $k$ (including or excluding 0 depending on the parity of $N$) by
order-preserving maps,
then under this identification the restriction of $\si$
becomes an element of \S. In fact it becomes an odd permutation because
the two
identification maps are odd. Thus we have described
a mapping $\si\ra{\cal F}(\si)$ from $\cO_{N+2}$ to $\cO_N$. The mapping
is
$2k+2$ to 1 and is clearly onto. It is also clear that $\ell_N({\cal
F}(\si))
\leq\ell_{N+2}(\si)$, from which it follows that $b_{N+2\;n}\leq
(2k+2)\,b_{N\,n}$
for all $n$. The assertion of the lemma follows upon using
(\ref{FN}).\sp

Theorem 1 tells us that
\begin{eqnarray*}
e^{2t^2}\varphi^e_n(2t^2)&=&\left\{
        \begin{array}{ll}
        D_{{n\over 2}}(t)^2, & \hspace{2.2em}n\ {\rm even,} \\&\\
        D_{{n-1\over 2}}(t) D_{{n+1\over 2}}(t), & \hspace{2.2em}n\ {\rm
odd,}
        \end{array}\right.\\
&&\\
e^{2t^2}\varphi^o_n(2t^2)&=&\left\{
        \begin{array}{ll}
        \left(1-U_{{n\over 2}}(t)^2\right) D_{{n\over 2}}(t)^2,
        &\> n\ {\rm even,} \\&\\
D_{{n-1\over 2}}(t) D_{{n+1\over 2}}(t), & \> n\ {\rm odd.}
        \end{array}\right.
\end{eqnarray*}
{}From the asymptotics (\ref{Dasym}) we find that if we set
\[{n\ov2}=2t+st^{1/3}\]
then
\[\lim_{t\ra\iy}\varphi^e_n(2t^2)=F(s)^2.\]
Setting $t=\sqrt{2k}/2$ gives
\[\lim_{k\ra\iy}\varphi^e_n(k)=F(s)^2.\]
Johansson's lemma tell us that when $N$ runs through the even integers
$2k$,
\[\lim_{N\ra\iy}F_N(2\sqrt{N}+2^{2/3}sN^{1/6})=\lim_{N\ra\iy}F_{N}(n)=
F(s)^2.\]

For $N$ running through the odd integers we obtain the same relation.
For this we
have to use the additional fact that $U_{n/2}(t)\ra 0$ as $t\ra\iy$. If
we
recall (\ref{UVD}) we see that this follows from the asymptotics of
$D_{n/2}$.
Thus the proof is complete.\sp

\begin{center}{\bf Acknowledgements}\end{center}
This work was supported in part by the National Science Foundation
through grants
DMS-9802122 (first author) and DMS-9732687 (second author). The authors
thank the
administration of the Mathematisches Forschungsinstitut Oberwolfach for
their
hospitality during the authors' visit under their Research in Pairs
program, when
the first results of the paper were obtained, and the
Volkswagen-Stiftung for its support of the program.\sp


\begin{thebibliography}{4}

\bibitem{BDJ} J.~Baik, P.~Deift and K.~Johansson, {\it On the
distribution
of the length of the longest increasing subsequence of random
permutations\/}, Preprint, math.CO/9810105.

\bibitem{DZ} P.~A.~Deift and X.~Zhou, {\it Asymptotics
for the Painlev\'e II equation\/}, Comm.\ Pure Appl.\ Math.\
{\bf 48} (1995) 277--337.

\bibitem{G} I.~M.~Gessel, {\it Symmetric functions and
P-recursiveness\/},
J.\ Combin. Th.\
Ser.~A, {\bf 53} (1990) 257--285.

\bibitem{H} M.~Hisakado, {\it Unitary matrix models and Painlev\'e
III\/},
Mod.\ Phys.\ Letts., {\bf A11} (1996) 3001--3010.

\bibitem{J} K.~Johansson, {\it The longest increasing subsequence
in a random permutation and a unitary random matrix model}, Math.\
Res.\ Letts., {\bf 5} (1998) 63--82.

\bibitem{Ni} F.~W.~Nijhoff and V.~G.~Papageorgiou,
{\it Similarity reductions of
integrable lattices and discrete analogues
of the Painlev\'e II equation\/}, Phys.\ Letts.\ {\bf 153A} (1991) 337--344.

\bibitem{Od} A.~M.~Odlyzko, B.~Poonen, H.~Widom and
H.~S.~Wilf, On the distribution of longest increasing
subsequences in random permutations, unpublished notes.

\bibitem{O} K.~Okamoto, {\it Polynomial Hamiltonians associated with
Painlev\'e equations},
Proc.\ Japan Acad., {\bf A56} (1980) 367--371.

\bibitem{Pe} V.~Periwal and D.~Shevitz, {\it Unitary-matrix models
as exactly solvable string theories\/}, Phys.\ Rev.\ Letts.\
{\bf 64} (1990) 1326--1329.

\bibitem{R} E.~M.~Rains, {\it Increasing subsequences and the classical
groups\/},
Elect.\ J.\ of Combinatorics, {\bf 5} \#R12 (1998).

\bibitem{TW} C.~A.~Tracy and H.~Widom, {\it Level-spacing distributions
and the Airy kernel\/},
 Comm.\ Math.\ Phys., {\bf 159} (1994) 151--174.

\end{thebibliography}
\end{document}